\documentclass[11pt]{article}
\usepackage{enumerate}
\usepackage{amssymb,a4wide,latexsym,makeidx,epsfig,fleqn}
\usepackage{amsthm}
\usepackage{amsmath}
\usepackage{enumerate}
\newtheorem{theorem}{Theorem}[section]

\newtheorem{lemma}[theorem]{Lemma}

\begin{document}
\textwidth 150mm \textheight 225mm
\title{Ordering of bicyclic graphs by matching energy
\thanks{ Supported by
the National Natural Science Foundation of China (Nos. 11171273 and 11601431) and the Seed Foundation of Innovation and Creation for Graduate Students in Northwestern Polytechnical University (No. Z2017190).}}
\author{{Xiangxiang Liu, Ligong Wang\footnote {Corresponding author.} and Peng Xiao}\\
{\small Department of Applied Mathematics, School of Science, Northwestern
Polytechnical University,}\\ {\small  Xi'an, Shaanxi 710072,
People's Republic
of China.}\\
{\small E-mail: xxliumath@163.com, lgwangmath@163.com, xiaopeng@sust.edu.cn
} }
\date{}
\maketitle
\begin{center}
\begin{minipage}{120mm}
\vskip 0.3cm
\begin{center}
{\small {\bf Abstract}}
\end{center}
{\small Let $G$ be a simple graph of order $n$ and $\mu_{1},\mu_{2},\ldots,\mu_{n}$ be the roots of its matching polynomial. The matching energy is defined as
the sum $\sum^{n}_{i=1}|\mu_{i}|$, which was introduced by Gutman and Wagner in 2012. In this paper, the graphs with the first five smallest matching energies among all bicyclic graphs for order $n>5$ are determined.

 }
\end{minipage}
\end{center}

\section{Introduction }
\label{sec:matching-introduction}

Let $G=(V,E)$ be a finite, connected, undirected and simple graph with vertex set $V=V(G)=\{v_{1},v_{2},\ldots,v_{n}\}$ and edge set $E=E(G)=\{e_{1},e_{2},\ldots,e_{m}\}$. A matching in a graph $G$ is a set of pairwise nonadjacent edges. A matching $M$ is called $k$-matching if the size of $M$ is $k$. Let $m(G,k)$ denote the number of $k$-matching of $G$, where $m(G,1)=m$, and $m(G,k)=0$ for $k>\lfloor\frac{n}{2}\rfloor$ or $k<0$. In addition, we define $m(G,0)=1$. The matching
polynomial of the graph $G$ is defined as
$$\alpha(G)=\alpha(G,x)=\sum\limits_{k\geq0}(-1)^{k}m(G,k)x^{n-2k}.$$

Let $\lambda_{1},\lambda_{2},\ldots,\lambda_{n}$ be the eigenvalues of the adjacency matrix of a graph $G$. The energy of graph $G$ \cite{GuI1} is defined as
$$E(G)=\sum\limits^{n}_{i=1}|\lambda_{i}|.$$

An important tool of graph energy is the Coulson integral formula \cite{GuI1} (with regard to $G$ being a tree $T$):
\begin{equation}\label{eq:matching-1}
E(T)=\frac{2}{\pi}\int^{\infty}_{0}\frac{1}{x^{2}}\ln\left[\sum\limits_{k\geq0}m(T,k)x^{2k}\right]dx.
\end{equation}

The energy of graphs has been widely studied by theoretical chemists and mathematicians. For details see the new book on graph energy \cite{GLI,LSG} and the reviews \cite{GuI2,GLZ}.

Recently, Gutman and Wagner \cite{GSI} defined the matching energy of a graph $G$. Let $G$ be a simple graph and
$\mu_{1},\mu_{2},\ldots,\mu_{n}$ be the zeros of its matching polynomial. Then
$$ME(G)=\sum^{n}_{i=1}|\mu_{i}|.$$

In view of Eq. \eqref{eq:matching-1}, the matching energy also has a beautiful formula as follows \cite{GSI}:
\begin{equation}\label{eq:matching-2}
ME(G)=\frac{2}{\pi}\int^{\infty}_{0}\frac{1}{x^{2}}\ln\left[\sum\limits_{k\geq0}m(G,k)x^{2k}\right]dx.
\end{equation}

By Eq. \eqref{eq:matching-2}, we know that the matching energy of a graph $G$ is a monotonically increasing function of $m(G,k)$. This means that if two graphs $G_{1}$ and $G_{2}$ satisfy $m(G_{1},k)\leq m(G_{2},k)$ for all $0\leq k\leq\lfloor\frac{n}{2}\rfloor$, then $ME(G_{1})\leq ME(G_{2})$. In addition, if $m(G_{1},k)<m(G_{2},k)$ for at least one $k$, then $ME(G_{1})<ME(G_{2})$.

We now introduce some elementary notations and terminologies that will be used in
the sequel. Let $G=(V,E)$ be a graph under our consideration. If $W\subseteq V(G)$, we denote by $G-W$ the subgraph of $G$ obtained by deleting the
vertices of $W$ and the edges incident with them. Similarly, if $E'\subseteq E(G)$, we denote by $G-E'$ the subgraph of $G$ obtained by deleting the edges of $E'$. If $W=\{v\}$ and $E'=\{xy\}$, we write $G-v$ and $G-xy$ instead of $G-\{v\}$ and $G-\{xy\}$, respectively.
We will use the notations $S_{n}$, $P_{n}$ and $C_{n}$ to denote the star, path and cycle on $n$ vertices, respectively. The union of the graphs $G_{1}$ and $G_{2}$, denoted by $G_{1}\cup G_{2}$, is the graph with vertex set $V(G_{1})\cup V(G_{2})$ and edge set $E(G_{1})\cup E(G_{2})$. Let $G_{1}$ and $G_{2}$ be two graphs with $V(G_{1})\cap V(G_{2})=\{v\}$. Let $G=G_{1}vG_{2}$ be a graph defined by $V(G)=V(G_{1})\cup V(G_{2})$ and $E(G)=E(G_{1})\cup E(G_{2})$. Let $u$ be a vertex of $G$, $N(u)$ or $N_{G}(u)$ denotes the neighborhood of $u$. We refer to Cvetkovi\'{c} et al. \cite{DMH} for undefined notations and terminologies.

\begin{figure}[htbp]
\begin{centering}
\includegraphics[scale=0.8]{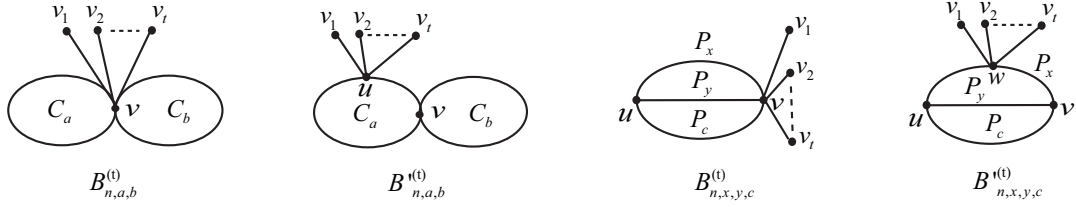}
\caption{Bicyclic graphs $B_{n,a,b}^{(t)}$, ${B'}_{n,a,b}^{(t)}$,
$B_{n,x,y,c}^{(t)}$ and ${B'}_{n,x,y,c}^{(t)}$.}\label{matching:Fig-1}
\end{centering}
\end{figure}

A bicyclic graph is a connected graph with $n$ vertices and $n+1$ edges. Let $\mathcal{B}(n)$ be the class of bicyclic graphs with $n$ vertices.
We shall use $B_{n,a,b}^{(t)}$ to denote the bicyclic graph constructed
by attaching $t$ pendent vertices to the vertex $v$ of $C_{a}vC_{b}$. Let ${B'}_{n,a,b}^{(t)}$ denote the bicyclic graph that is obtained by attaching $t$ pendent vertices to one vertex except $v$ of $C_{a}vC_{b}$. Three internal disjoint paths $P_{x}$, $P_{y}$ and $P_{c}$ possessing common end vertices $u$, $v$ form a bicyclic graph denoted by $B_{x,y,c}$. We shall use $B_{n,x,y,c}^{(t)}$ to denote the bicyclic graph constructed by attaching $t$ pendent vertices to the vertex $v$ of $B_{x,y,c}$. Let ${B'}_{n,x,y,c}^{(t)}$ denote the bicyclic graph that is obtained by attaching $t$ pendent vertices to one vertex except $v$ and $u$ of $B_{x,y,c}$. Bicyclic graphs $B_{n,a,b}^{(t)}$, ${B'}_{n,a,b}^{(t)}$,
$B_{n,x,y,c}^{(t)}$ and ${B'}_{n,x,y,c}^{(t)}$ are depicted in Figure \ref{matching:Fig-1}.

The study on extremal matching energy is interesting. In \cite{GSI}, Gutman and Wagner gave some elementary results on the matching energy and obtained the unicyclic graphs with the minimal and maximal matching energy. For the bicyclic graphs, Ji et al. \cite{JLS} obtained the graphs with the minimal and maximal matching energy.
In \cite{JMS}, Ji and Ma obtained tricyclic graph with maximum matching energy. In \cite{ZLL}, Zou and Li characterized the bicyclic graph with given girth having minimum matching energy. For more results about matching energy, see \cite{CLS,CLS1,LYS,WSW,XDZ}.

This paper is organized as follows: In Section 2, we give some preliminary results, which will be used in the following discussion. For $n>5$, the graphs with the first five smallest matching energies among all bicyclic graphs of order $n$ are determined in Section 3.

\section{Preliminary Results}
\label{sec:matching-Preliminary}

In this section, we shall give some elementary results which will be used in the following.

\begin{figure}[htbp]
\begin{centering}
\includegraphics[scale=0.6]{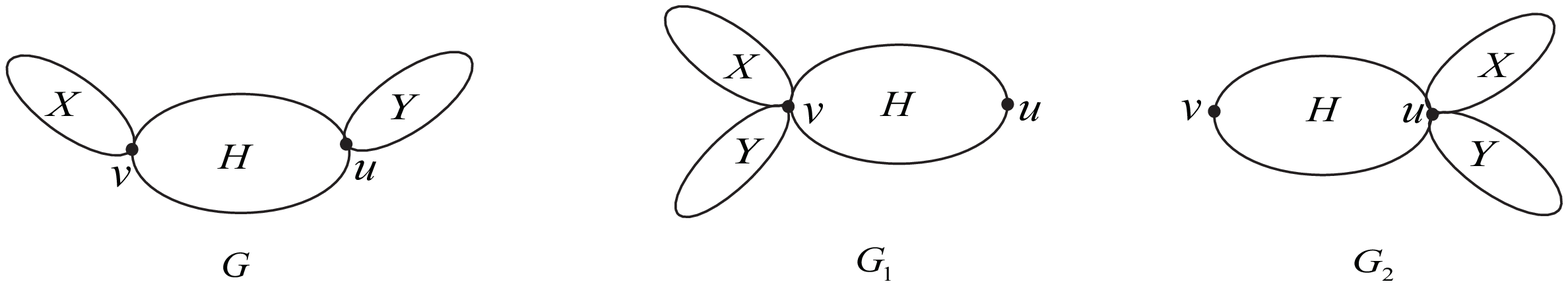}
\caption{Graphs $G$, $G_{1}$ and $G_{2}$.}\label{matching:Fig-2}
\end{centering}
\end{figure}

\noindent\begin{lemma}\label{le:matching-1} (\cite{FEJ,GuI3})
Let $G=(V,E)$ be a graph.

(i) If $uv\in E(G)$, then $m(G,k)=m(G-uv,k)+m(G-u-v,k-1)$;

(ii) If $u\in V(G)$, then $m(G,k)=m(G-u,k)+\sum\limits_{v\in N(u)} m(G-u-v,k-1)$.
\end{lemma}
\noindent\begin{lemma}\label{le:matching-2} (\cite{GuI})
If $m(G_{1},k)\geq m(G_{2},k)$, then $m(G_{1}\cup H,k)\geq m(G_{2}\cup H,k)$, where $H$ is an arbitrary graph.
\end{lemma}
\noindent\begin{lemma}\label{le:matching-3} (\cite{CLS})
Let $G$ be a simple graph and $H$ be a subgraph (resp. proper subgraph) of $G$. Then $m(G,k)\geq m(H,k)$ (resp. $m(G,k)>m(H,k)$).
\end{lemma}
\noindent\begin{lemma}\label{le:matching-4} (\cite{GSI})
Suppose that $G$ is a connected graph and $T$ is an induced subgraph of $G$ such that $T$ is a tree and $T$ is connected to the rest of $G$ only by a cut vertex $v$. If $T$ is replaced by a star of the same order, centered at $v$, then the matching energy decreases (unless $T$ is already such a star). If $T$ is replaced by a path, with one end at $v$, then the matching energy increases (unless $T$ is already such a path).
\end{lemma}
\noindent\begin{lemma}\label{le:matching-5}
Let $H$, $X$, $Y$ be three connected graphs disjoint in pair. Suppose that $u$, $v$ are two vertices of $H$, $v'$ is a vertex of $X$, $u'$ is a vertex of $Y$. Let $G$ be the graph obtained from $H$, $X$, $Y$ by identifying $v$ with $v'$ and $u$ with $u'$, respectively. Let $G_{1}$ be the graph obtained
from $H$, $X$, $Y$ by identifying vertices $v$, $v'$, $u'$ and $G_{2}$ be the graph obtained from $H$, $X$, $Y$ by identifying
vertices $u$, $v'$, $u'$; see Figure \ref{matching:Fig-2}. Then
$ME(G_{1})<ME(G)$ and $ME(G_{2})<ME(G).$
\end{lemma}
\noindent {\bf Proof.} By Lemma \ref{le:matching-1}, we have
\begin{align*}
 m(G,k)=&m(G-v,k)+\sum\limits_{w\in N(v)}m(G-v-w,k-1),
\end{align*}
\begin{align*}
 m(G_{1},k)=&m(G_{1}-v,k)+\sum\limits_{w'\in N(v)}m(G_{1}-v-w',k-1).
\end{align*}
Since $G_{1}-v$ is a subgraph of $G-v$ and $G_{1}-v-w'$ is a subgraph of $G-v-w$.
By Lemma \ref{le:matching-3}, we have $m(G_{1}-v,k)<m(G-v,k)$ and $m(G_{1}-v-w',k-1)<m(G-v-w,k-1)$. Hence, $ME(G_{1})<ME(G)$.

Similarly, we have $ME(G_{2})<ME(G).$   \hfill$\blacksquare$

\section{Main Results}
\label{sec:matching-Main Results}

\noindent\begin{lemma}\label{th:matching-1}
For positive integers $a$, $b$ and $t$ with $a,b\geq3$, $ME(B_{n,a,b}^{(t)})\leq ME({B'}_{n,a,b}^{(t)})$.
\end{lemma}
\noindent {\bf Proof.} Let $B_{1}=B_{n-t,a,b}^{(0)}$. By Lemma \ref{le:matching-1}, we have
\begin{align*}
 m(B_{n,a,b}^{(t)},k)=&m(B_{n,a,b}^{(t)}-v_{1},k)+m(B_{n,a,b}^{(t)}-v_{1}-v,k-1)\\
                   =&m(B_{n-1,a,b}^{(t-1)},k)+m(P_{a-1}\cup P_{b-1},k-1)\\
                   =&m(B_{n-2,a,b}^{(t-2)},k)+2m(P_{a-1}\cup P_{b-1},k-1)\\
                   =&\cdots \\
                   =&m(B_{n-t,a,b}^{(0)},k)+tm(P_{a-1}\cup P_{b-1},k-1)\\
                   =&m(B_{n-t,a,b}^{(0)}-v,k)+\sum\limits_{w\in N_{B_{1}}(v)} m(B_{n-t,a,b}^{(0)}-v-w,k-1)+tm(P_{a-1}\cup P_{b-1},k-1)\\
                   =&m(P_{a-1}\cup P_{b-1},k)+2m(P_{a-2}\cup P_{b-1},k-1)\\&+2m(P_{a-1}\cup P_{b-2},k-1)+tm(P_{a-1}\cup P_{b-1},k-1).
\end{align*}

Let $Q={B'}_{n,a,b}^{(t)}-\{u,v_{1},\ldots,v_{t}\}$. For two vertices
$u,\ v\in V(C_{a})$, $P_{x}$ denotes the subpath of $C_{a}$ from $u$ to $v$ and $P_{y}$ denotes the subpath from $v$ to $u$ in the reversed direction of $C_{a}$ and $a=x+y-2$. Let $B_{2}={B'}_{n-t,a,b}^{(0)}$. By Lemma \ref{le:matching-1}, we have
\begin{align*}
 m({B'}_{n,a,b}^{(t)},k)=&m({B'}_{n,a,b}^{(t)}-v_{1},k)+m({B'}_{n,a,b}^{(t)}-v_{1}-u,k-1)\\
                   =&m({B'}_{n-1,a,b}^{(t-1)},k)+m(Q,k-1)\\
                   =&m({B'}_{n-2,a,b}^{(t-2)},k)+m(Q,k-1)\\
                   =&\cdots \\
                   =&m({B'}_{n-t,a,b}^{(0)},k)+tm(Q,k-1)\\
                    =&m({B'}_{n-t,a,b}^{(0)}-v,k)+\sum\limits_{r\in N_{B_{2}}(v)} m({B'}_{n-t,a,b}^{(0)}-v-r,k-1)+tm(Q,k-1)\\
 \end{align*}
 \begin{align*}
                   =&m(P_{a-1}\cup P_{b-1},k)+2m(P_{a-2}\cup P_{b-1},k-1)\\&+2m(P_{a-1}\cup P_{b-2},k-1)+tm(Q,k-1),
\end{align*}
where
\begin{align*}
 m(Q,k-1)=&m(Q-v,k-1)+\sum\limits_{w'\in N_{Q}(v)}m(Q-v-w',k-2)\\
         =&m(P_{x-2}\cup P_{y-2}\cup P_{b-1},k-1)+m(P_{x-3}\cup P_{y-2}\cup P_{b-1},k-2)\\
         &+m(P_{x-2}\cup P_{y-3}\cup P_{b-1},k-2)+2m(P_{x-2}\cup P_{y-2}\cup P_{b-2},k-2)\\
         =&m(P_{x-1}\cup P_{y-2}\cup P_{b-1},k-1)+m(P_{x-2}\cup P_{y-3}\cup P_{b-1},k-2)\\
         &+2m(P_{x-2}\cup P_{y-2}\cup P_{b-2},k-2).
\end{align*}

Note that $a=x+y-2$, by Lemma \ref{le:matching-1} we have
\begin{align*}
 m(P_{a-1}\cup P_{b-1},k-1)=&m(P_{x-1}\cup P_{y-2}\cup P_{b-1},k-1)+m(P_{x-2}\cup P_{y-3}\cup P_{b-1},k-2).
\end{align*}

Combining with above results,
$$m({B'}_{n,a,b}^{(t)},k)-m(B_{n,a,b}^{(t)},k)=2tm(P_{x-2}\cup P_{y-2}\cup P_{b-2},k-2)\geq0.$$
Hence, $ME(B_{n,a,b}^{(t)})\leq ME({B'}_{n,a,b}^{(t)})$.  \hfill$\blacksquare$

\begin{figure}[htbp]
\begin{centering}
\includegraphics[scale=0.85]{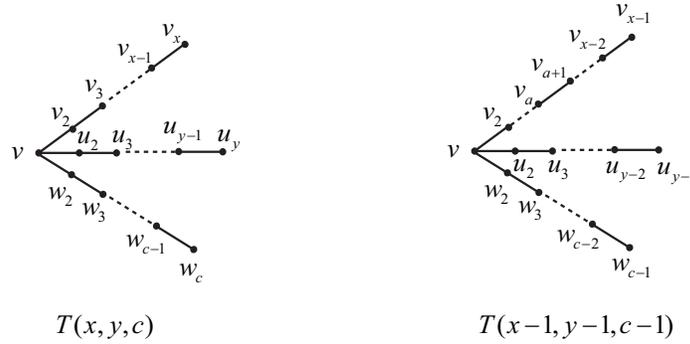}
\caption{Trees $T(x,y,c)$ and $T(x-1,y-1,c-1)$.}\label{matching:Fig-3}
\end{centering}
\end{figure}

Denoted by $T(x,y,c)$ the tree with exactly one vertex $v$ of degree 3, and $T(x,y,c)-v=P(x-1)\cup P(y-1)\cup P(c-1)$; see Figure \ref{matching:Fig-3}.

\noindent\begin{lemma}\label{th:matching-2}
For positive integers $x$, $y$, $c$ and $t$ with $x\geq 3,\ yc\geq6$, $ME(B_{n,x,y,c}^{(t)})\leq ME({B'}_{n,x,y,c}^{(t)})$.
\end{lemma}
\noindent {\bf Proof.} Let $T=B_{n,x,y,c}^{(t)}-\{v,v_{1},\ldots,v_{t}\}=T(x-1,y-1,c-1)$ and $B_{3}=B_{n-t,x,y,c}^{(0)}$. By Lemma \ref{le:matching-1}, we have
\begin{align*}
 m(B_{n,x,y,c}^{(t)},k)=&m(B_{n,x,y,c}^{(t)}-v_{1},k)+m(B_{n,x,y,c}^{(t)}-v_{1}-v,k-1)\\
                   =&m(B_{n-1,x,y,c}^{(t-1)},k)+m(T,k-1)\\
                   =&m(B_{n-2,x,y,c}^{(t-2)},k)+2m(T,k-1)\\
                   =&\cdots \\
                   =&m(B_{n-t,x,y,c}^{(0)},k)+tm(T,k-1)\\
 \end{align*}
 \begin{align*}
                   =&m(B_{n-t,x,y,c}^{(0)}-v,k)+\sum\limits_{r\in N_{B_{3}}(v)} m(B_{n-t,x,y,c}^{(0)}-v-r,k-1)+tm(T,k-1)\\
                   =&m(T,k)+\sum\limits_{r\in N_{B_{3}}(v)} m(T-r,k-1)+tm(T,k-1).
\end{align*}
where
\begin{align*}
 m(T,k)=&m(T-u,k)+\sum\limits_{r'\in N_{T}(u)}m(T-u-r',k-1)\\
         =&m(P_{x-2}\cup P_{y-2}\cup P_{c-2},k)+m(P_{x-3}\cup P_{y-2}\cup P_{c-2},k-1)\\
         &+m(P_{x-2}\cup P_{y-3}\cup P_{c-2},k-1)+m(P_{x-2}\cup P_{y-2}\cup P_{c-3},k-1)\\
         =&m(P_{x-1}\cup P_{y-2}\cup P_{c-2},k)+m(P_{x-2}\cup P_{y-3}\cup P_{c-2},k-1)\\
         &+m(P_{x-2}\cup P_{y-2}\cup P_{c-3},k-1).
\end{align*}
Furthermore, let $N_{B_{3}}(v)=\{r_{1},r_{2},r_{3}\}$. By Lemma \ref{le:matching-1}, we have
\begin{align*}
 m(T-r_{1},k-1)=&m(T-r_{1}-u,k-1)+\sum\limits_{r'\in N_{T-r_{1}}(u)}m(T-r_{1}-u-r',k-2)\\
         =&m(P_{x-3}\cup P_{y-2}\cup P_{c-2},k-1)+m(P_{x-4}\cup P_{y-2}\cup P_{c-2},k-2)\\
         &+m(P_{x-3}\cup P_{y-3}\cup P_{c-2},k-2)+m(P_{x-3}\cup P_{y-2}\cup P_{c-3},k-2)\\
         =&m(P_{x-2}\cup P_{y-2}\cup P_{c-2},k-1)+m(P_{x-3}\cup P_{y-2}\cup P_{c-3},k-2)\\
         &+m(P_{x-3}\cup P_{y-3}\cup P_{c-2},k-2).
\end{align*}
Similarly, we have
\begin{align*}
 m(T-r_{2},k-1)=&m(P_{x-2}\cup P_{y-2}\cup P_{c-2},k-1)+m(P_{x-3}\cup P_{y-2}\cup P_{c-3},k-2)\\&+m(P_{x-2}\cup P_{y-3}\cup P_{c-3},k-2).
\end{align*}
\begin{align*}
 m(T-r_{3},k-1)=&m(P_{x-2}\cup P_{y-2}\cup P_{c-2},k-1)+m(P_{x-3}\cup P_{y-3}\cup P_{c-2},k-2)\\&+m(P_{x-2}\cup P_{y-3}\cup P_{c-3},k-2).
\end{align*}

Combining with above results,
\begin{align*}
 m(B_{n,x,y,c}^{(t)},k)=&m(P_{x-1}\cup P_{y-2}\cup P_{c-2},k)+m(P_{x-2}\cup P_{y-3}\cup P_{c-2},k-1)\\
                  &+m(P_{x-2}\cup P_{y-2}\cup P_{c-3},k-1)+3m(P_{x-2}\cup P_{y-2}\cup P_{c-2},k-1)\\
                  &+2m(P_{x-3}\cup P_{y-2}\cup P_{c-3},k-2)+2m(P_{x-3}\cup P_{y-3}\cup P_{c-2},k-2)\\
                  &+2m(P_{x-2}\cup P_{y-3}\cup P_{c-3},k-2)+tm(P_{x-1}\cup P_{y-2}\cup P_{c-2},k-1)\\
                  &+tm(P_{x-2}\cup P_{y-3}\cup P_{c-2},k-2)+tm(P_{x-2}\cup P_{y-2}\cup P_{c-3},k-2).
\end{align*}

Let $H={B'}_{n,x,y,c}^{(t)}-\{w,v_{1},\ldots,v_{t}\}$. Similarly, we have
\begin{align*}
 m({B'}_{n,x,y,c}^{(t)},k)=&m(P_{x-1}\cup P_{y-2}\cup P_{c-2},k)+m(P_{x-2}\cup P_{y-3}\cup P_{c-2},k-1)\\
                  &+m(P_{x-2}\cup P_{y-2}\cup P_{c-3},k-1)+3m(P_{x-2}\cup P_{y-2}\cup P_{c-2},k-1)\\
                  &+2m(P_{x-3}\cup P_{y-2}\cup P_{c-3},k-2)+2m(P_{x-3}\cup P_{y-3}\cup P_{c-2},k-2)\\
                  &+2m(P_{x-2}\cup P_{y-3}\cup P_{c-3},k-2)+tm(H,k-1).
\end{align*}

Therefore, $m({B'}_{n,x,y,c}^{(t)},k)- m(B_{n,x,y,c}^{(t)},k)=tm(H,k-1)-tm(T,k-1).$ Let $H-v=T(a,y-1,c-1)\cup P_{x-a-2}$.

By lemma \ref{le:matching-1}, we have
\begin{align*}
 m(H,k-1)=&m(H-v,k-1)+\sum\limits_{w'\in N_{H}(v)}m(H-v-w',k-2)\\
         =&m(T(a,y-1,c-1)\cup P_{x-a-2},k-1)+m(T(a,y-1,c-1)\cup P_{x-a-3},k-2)\\
         &+m(T(a,y-2,c-1)\cup P_{x-a-2},k-2)+m(T(a,y-1,c-2)\cup P_{x-a-2},k-2)\\
         =&m(T(a,y-1,c-1)\cup P_{x-a-1},k-1)+m(T(a,y-2,c-1)\cup P_{x-a-2},k-2)\\
         &+m(T(a,y-1,c-2)\cup P_{x-a-2},k-2).
\end{align*}
Let $e=v_{a}v_{a+1}\in T$; see Figure \ref{matching:Fig-3}. Then
\begin{align*}
 m(T,k-1)=&m(T-e,k-1)+m(T-u'-v',k-2)\\
         =&m(T(a,y-1,c-1)\cup P_{x-a-1},k-1)+m(T(a-1,y-1,c-1)\cup P_{x-a-2},k-2).
\end{align*}
Hence,
 \begin{align*}
 &m(H,k-1)-m(T,k-1)\\
 =&m(T(a,y-2,c-1)\cup P_{x-a-2},k-2)+m(T(a,y-1,c-2)\cup P_{x-a-2},k-2)\\
 &-m(T(a-1,y-1,c-1)\cup P_{x-a-2},k-2)\\
 =&m(P_{c-2}\cup P_{a+y-3}\cup P_{x-a-2},k-2)+m(P_{c-3}\cup P_{a-1}\cup P_{y-3}\cup P_{x-a-2},k-3)\\
 &+m(T(a,y-1,c-2)\cup P_{x-a-2},k-2)\\
  &-m(P_{c-2}\cup P_{a+y-3}\cup P_{x-a-2},k-2)-m(P_{c-3}\cup P_{a-2}\cup P_{y-2}\cup P_{x-a-2},k-3)\\
 =&m(P_{c-3}\cup P_{a-1}\cup P_{y-3}\cup P_{x-a-2},k-3)+m(P_{c-3}\cup P_{a+y-2} \cup P_{x-a-2},k-2)\\
 &+m(P_{c-4}\cup P_{a-1}\cup P_{y-2}\cup P_{x-a-2},k-3)-m(P_{c-3}\cup P_{a-2}\cup P_{y-2}\cup P_{x-a-2},k-3)\\
 =&m(P_{c-3}\cup P_{a-1}\cup P_{y-3}\cup P_{x-a-2},k-3)+m(P_{c-3}\cup P_{a-1}\cup P_{y-1} \cup P_{x-a-2},k-2)\\
 &+m(P_{c-4}\cup P_{a-1}\cup P_{y-2}\cup P_{x-a-2},k-3).
\end{align*}
By Lemmas \ref{le:matching-2} and \ref{le:matching-3}, we have $m(H,k-1)-m(T,k-1)\geq0.$ That is to say ${B'}_{n,x,y,c}^{(t)},k)- m(B_{n,x,y,c}^{(t)},k)\geq0$.

Thus, for positive integers $x$, $y$, $c$ and $t$ with $x\geq 3,\ yc\geq6$, $ME(B_{n,x,y,c}^{(t)})\leq ME({B'}_{n,x,y,c}^{(t)})$.

 \hfill$\blacksquare$

\noindent\begin{lemma}\label{th:matching-3}
Let $G\in \mathcal{B}(n)$. Then

(i) If $G$ contains exactly two cycles, say $C_{a}$ and $C_{b}$, then $ME(G)\geq ME(B_{n,a,b}^{(t)})$, the equality holds if and only if $G\cong B_{n,a,b}^{(t)}$.

(ii) If $G$ contains exactly three cycles, say $B_{x,y,c}$, then $ME(G)\geq ME(B_{n,x,y,c}^{(t)})$, the equality holds if and only if $G\cong B_{n,x,y,c}^{(t)}$.
\end{lemma}
\noindent {\bf Proof.} (i) If $G$ contains exactly two cycles, say $C_{a}$ and $C_{b}$. Assume that $C_{a}$ connects $C_{b}$ by a path $P_{l+2}$ with $l\geq -1$. It is easy to know that $C_{a}$ and $C_{b}$ have exactly one vertex in common  when $l=-1$. Let $C_{a}=u_{1}u_{2}\ldots u_{a}u_{1}$, $C_{b}=v_{1}v_{2}\ldots v_{b}v_{1}$ and $P_{l+2}=u_{1}w_{1}w_{2}\ldots w_{l}v_{1}$. Set $V^{\ast}(G)=\{u_{i}: d(u_{i})\geq3,2\leq i\leq a\}\cup \{v_{i}: d(v_{i})\geq3,2\leq i\leq b\}\cup\{w_{i}: d(w_{i})\geq3,1\leq i\leq l\}\cup \{u_{1}: d(u_{1})\geq4\}\cup\{v_{1}: d(v_{1})\geq4\}$. Suppose that $|V^{\ast}(G)|=k$, and relabel the vertices in $V^{\ast}(G)$ as $\{r_{1},r_{2},\ldots,r_{k}\}$.

Let $r_{i}\in V^{\ast}(G)$ and $T_{i}$ be a subtree of $G-E(C_{a}\cup C_{b}\cup P_{l+2})$ which contains $r_{i}$ and $|V(T_{i})|=p_{i}+1$. Let $V'(T_{i})=V(T_{i})-r_{i}$. Denote $H=G-V'(T_{i})$.

Then $G=Hr_{i}T_{i}$. By Lemma \ref{le:matching-4}, we have $ME(Hr_{i}T_{i})\geq ME(Hr_{i}S_{p_{i}+1})$. Thus repeatedly using Lemma \ref{le:matching-4}, we have
$ME(G)\geq ME(B_{n}(p_{1},p_{2},\ldots,p_{k}))$, where $B_{n}(p_{1},p_{2},\ldots,p_{k})$ is a bicyclic graph with $n$ vertices created from $C_{a}u_{1}P_{l+2}v_{1}C_{b}$ by attaching $p_{i}$ pendent vertices to $r_{i}\in V^{\ast}(G)$, $1\leq i\leq k$, respectively.

Let $X=S_{p_{i}+1}$, $Y=S_{p_{j}+1}$, $V'(X)=V(X)-r_{i}$ and $V'(Y)=V(Y)-r_{j}$. Denote $H'=G-V'(X)-V'(Y)$. Then $B_{n}(p_{1},p_{2},\ldots,p_{k})=Xr_{i}H'r_{j}Y$, where $r_{i}$ and $r_{j}$ are centers of $X$ and $Y$, respectively. By Lemma \ref{le:matching-5}, we have
\begin{align*}
 ME(G)\geq ME(B_{n}(p_{1},\ldots,p_{i},\ldots,p_{j},\ldots,p_{k}))>ME(B_{n}(p_{1},\ldots,p_{i+j},\ldots,0,\ldots,p_{k})),
\end{align*}
or
\begin{align*}
ME(G)\geq ME(B_{n}(p_{1},\ldots,p_{i},\ldots,p_{j},\ldots,p_{k}))>ME(B_{n}(p_{1},\ldots,0,\ldots,p_{i+j},\ldots,p_{k})).
\end{align*}
Repeatedly using above step, and by Lemmas \ref{le:matching-1}, \ref{le:matching-2} and \ref{le:matching-3}, we obtain either
\begin{align*}
ME(G)\geq ME(B_{n}(p_{1},\ldots,p_{i},\ldots,p_{j},\ldots,p_{k}))>\cdots>ME(B_{n,a,b}^{(t)}),
\end{align*}
or
\begin{align*}
ME(G)\geq
ME(B_{n}(p_{1},\ldots,p_{i},\ldots,p_{j},\ldots,p_{k}))>\cdots>ME({B'}_{n,a,b}^{(t)}).
\end{align*}

Combining with Lemma \ref{th:matching-1}, if $G$ contains exactly two cycles, then we have
$ME(G)\geq ME(B_{n,a,b}^{(t)})$, the equality holds if and only if $G\cong B_{n,a,b}^{(t)}$.

(ii) If $G$ contains exactly three cycles, say $B_{x,y,c}$. Let $P_{x}=u_{1}u_{2}\ldots u_{x}$, $P_{y}=u_{1}v_{2}\ldots v_{y-1}u_{x}$ and $P_{c}=u_{1}w_{2}\ldots w_{c-1}u_{x}$. Set $V^{\ast}(G)=\{u_{i}: d(u_{i})\geq3,2\leq i\leq x-1\}\cup \{v_{i}: d(v_{i})\geq3,2\leq i\leq y-1\}\cup\{w_{i}: d(w_{i})\geq3,2\leq i\leq c-1\}\cup \{u_{1}: d(u_{1})\geq4\}\cup\{u_{x}: d(u_{x})\geq4\}$. Assume that $|V^{\ast}(G)|=k$, and relabel the vertices in $V^{\ast}(G)$ as $\{r_{1},r_{2},\ldots,r_{k}\}$.

Let $r_{i}\in V^{\ast}(G)$ and $T_{i}$ be a subtree of $G-E(B_{x,y,c})$ which contains $r_{i}$ and $|V(T_{i})|=p_{i}+1$. Let $V'(T_{i})=V(T_{i})-r_{i}$. Denote $H=G-V'(T_{i})$.

Then $G=Hr_{i}T_{i}$. By Lemma \ref{le:matching-4}, we have $ME(Hr_{i}T_{i})\geq ME(Hr_{i}S_{p_{i}+1})$. Thus repeated using Lemma \ref{le:matching-4}, we have
$ME(G)\geq ME(B_{n}'(p_{1},p_{2},\ldots,p_{k}))$, where $B_{n}'(p_{1},p_{2},\ldots,p_{k})$ is a bicyclic graph with $n$ vertices created from $B_{x,y,c}$ by attaching $p_{i}$ pendent vertices to $r_{i}\in V^{\ast}(G)$, $1\leq i\leq k$, respectively.

Let $X=S_{p_{i}+1}$, $Y=S_{p_{j}+1}$, $V'(X)=V(X)-r_{i}$ and $V'(Y)=V(Y)-r_{j}$. Denote $H'=G-V'(X)-V'(Y)$. Then $B_{n}'(p_{1},p_{2},\ldots,p_{k})=Xr_{i}H'r_{j}Y$, where $r_{i}$ and $r_{j}$ are centers of $X$ and $Y$, respectively. By Lemma \ref{le:matching-5}, we have
\begin{align*}
ME(G)\geq ME(B_{n}'(p_{1},\ldots,p_{i},\ldots,p_{j},\ldots,p_{k}))>ME(B_{n}'(p_{1},\ldots,p_{i+j},\ldots,0,\ldots,p_{k})),
\end{align*}
or
\begin{align*}
ME(G)\geq ME(B_{n}'(p_{1},\ldots,p_{i},\ldots,p_{j},\ldots,p_{k}))>ME(B_{n}'(p_{1},\ldots,0,\ldots,p_{i+j},\ldots,p_{k})).
\end{align*}
Repeatedly using above step, we obtain either
\begin{align*}
ME(G)\geq ME(B_{n}'(p_{1},\ldots,p_{i},\ldots,p_{j},\ldots,p_{k}))>\cdots>ME(B_{n,x,y,c}^{(t)}),
\end{align*}
or
\begin{align*}
ME(G)\geq
ME(B_{n}'(p_{1},\ldots,p_{i},\ldots,p_{j},\ldots,p_{k}))>\cdots>ME({B'}_{n,x,y,c}^{(t)}).
\end{align*}
Combining with Lemma \ref{th:matching-2}, if $G$ contains exactly three cycles, then we have
$ME(G)\geq ME(B_{n,x,y,c}^{(t)})$, the equality holds if and only if $G\cong B_{n,x,y,c}^{(t)}$. \hfill$\blacksquare$

\noindent\begin{theorem}\label{th:matching-4}
For positive integers $a$, $b$ and $t$ with $a\geq 4,b\geq3$, $ME(B_{n,a-1,b}^{(t+1)})<ME(B_{n,a,b}^{(t)})$.
\end{theorem}
\noindent {\bf Proof.} By the proof of Lemma \ref{th:matching-1}, we have
\begin{align*}
 m(B_{n,a-1,b}^{(t+1)},k)-m(B_{n,a,b}^{(t)},k)
                   =&m(P_{a-2}\cup P_{b-1},k)+2m(P_{a-3}\cup P_{b-1},k-1)\\
                   &+2m(P_{a-2}\cup P_{b-2},k-1)+(t+1)m(P_{a-2}\cup P_{b-1},k-1)\\
                   &-m(P_{a-1}\cup P_{b-1},k)-2m(P_{a-2}\cup P_{b-1},k-1)\\
                   &-2m(P_{a-1}\cup P_{b-2},k-1)-tm(P_{a-1}\cup P_{b-1},k-1)\\
                   =&m(P_{a-2}\cup P_{b-1},k)+2m(P_{a-3}\cup P_{b-1},k-1)\\
                   &+2m(P_{a-2}\cup P_{b-2},k-1)+(t+1)m(P_{a-2}\cup P_{b-1},k-1)\\
                   &-m(P_{a-2}\cup P_{b-1},k)-m(P_{a-2}\cup P_{b-1},k-1)\\
                   &-m(P_{a-3}\cup P_{b-1},k-1)-m(P_{a-2}\cup P_{b-1},k-1)\\
                   &-2m(P_{a-1}\cup P_{b-2},k-1)-tm(P_{a-1}\cup P_{b-1},k-1)\\
                   =&m(P_{a-3}\cup P_{b-1},k-1)-m(P_{a-2}\cup P_{b-1},k-1)\\
                   &+2m(P_{a-2}\cup P_{b-2},k-1)-2m(P_{a-1}\cup P_{b-2},k-1)\\
                   &+tm(P_{a-2}\cup P_{b-1},k-1)-tm(P_{a-1}\cup P_{b-1},k-1).
\end{align*}
By Lemmas \ref{le:matching-2} and \ref{le:matching-3}, we have $m(B_{n,a-1,b}^{(t+1)},k)-m(B_{n,a,b}^{(t)},k)<0$. So $ME(B_{n,a-1,b}^{(t+1)})<ME(B_{n,a,b}^{(t)})$. \hfill$\blacksquare$

\noindent\begin{theorem}\label{th:matching-5}
For positive integers $x$, $y$, $c$ and $t$ with $x\geq 4,\ y,c\geq2,\ yc\geq6$, $ME(B_{n,x-1,y,c}^{(t+1)})$ $<ME(B_{n,x,y,c}^{(t)}).$
\end{theorem}
\noindent {\bf Proof.} By the proof of Lemma \ref{th:matching-2}, we have
\begin{align*}
 &m(B_{n,x-1,y,c}^{(t+1)},k)-m(B_{n,x,y,c}^{(t)},k)\\
                   =&m(P_{x-2}\cup P_{y-2}\cup P_{c-2},k)+m(P_{x-3}\cup P_{y-3}\cup P_{c-2},k-1)\\
                  &+m(P_{x-3}\cup P_{y-2}\cup P_{c-3},k-1)+3m(P_{x-3}\cup P_{y-2}\cup P_{c-2},k-1)\\
                  &+2m(P_{x-4}\cup P_{y-2}\cup P_{c-3},k-2)+2m(P_{x-4}\cup P_{y-3}\cup P_{c-2},k-2)\\
                  &+2m(P_{x-3}\cup P_{y-3}\cup P_{c-3},k-2)+(t+1)m(P_{x-2}\cup P_{y-2}\cup P_{c-2},k-1)\\
                  &+(t+1)m(P_{x-3}\cup P_{y-3}\cup P_{c-2},k-2)+(t+1)m(P_{x-3}\cup P_{y-2}\cup P_{c-3},k-2)\\
                   &-m(P_{x-1}\cup P_{y-2}\cup P_{c-2},k)-3m(P_{x-2}\cup P_{y-2}\cup P_{c-2},k-1)\\
                   &-m(P_{x-2}\cup P_{y-3}\cup P_{c-2},k-1)-2m(P_{x-3}\cup P_{y-3}\cup P_{c-2},k-2)\\
                  &-m(P_{x-2}\cup P_{y-2}\cup P_{c-3},k-1)-2m(P_{x-3}\cup P_{y-2}\cup P_{c-3},k-2)\\
                  &-2m(P_{x-2}\cup P_{y-3}\cup P_{c-3},k-2)-tm(P_{x-1}\cup P_{y-2}\cup P_{c-2},k-1)\\
                  &-tm(P_{x-2}\cup P_{y-3}\cup P_{c-2},k-2)-tm(P_{x-2}\cup P_{y-2}\cup P_{c-3},k-2)\\
   \end{align*}
 \begin{align*}
                  =&m(P_{x-2}\cup P_{y-2}\cup P_{c-2},k)+m(P_{x-3}\cup P_{y-3}\cup P_{c-2},k-1)\\
                  &+m(P_{x-3}\cup P_{y-2}\cup P_{c-3},k-1)+3m(P_{x-3}\cup P_{y-2}\cup P_{c-2},k-1)\\
                  &+2m(P_{x-4}\cup P_{y-2}\cup P_{c-3},k-2)+2m(P_{x-4}\cup P_{y-3}\cup P_{c-2},k-2)\\
                  &+2m(P_{x-3}\cup P_{y-3}\cup P_{c-3},k-2)+(t+1)m(P_{x-2}\cup P_{y-2}\cup P_{c-2},k-1)\\
                  &+(t+1)m(P_{x-3}\cup P_{y-3}\cup P_{c-2},k-2)+(t+1)m(P_{x-3}\cup P_{y-2}\cup P_{c-3},k-2)\\
                  &-[m(P_{x-2}\cup P_{y-2}\cup P_{c-2},k)+m(P_{x-2}\cup P_{y-2}\cup P_{c-2},k-1)\\
                  &+2m(P_{x-2}\cup P_{y-2}\cup P_{c-2},k-1)+m(P_{x-3}\cup P_{y-2}\cup P_{c-2},k-1)]\\
                   &-[m(P_{x-3}\cup P_{y-3}\cup P_{c-2},k-1)+m(P_{x-3}\cup P_{y-3}\cup P_{c-2},k-2)\\
                  &+m(P_{x-3}\cup P_{y-3}\cup P_{c-2},k-2)+m(P_{x-4}\cup P_{y-3}\cup P_{c-2},k-2)]\\
                  &-[m(P_{x-3}\cup P_{y-2}\cup P_{c-3},k-1)+m(P_{x-3}\cup P_{y-2}\cup P_{c-3},k-2)\\
                  &+m(P_{x-3}\cup P_{y-2}\cup P_{c-3},k-2)
                  +m(P_{x-4}\cup P_{y-2}\cup P_{c-3},k-2)]\\
                  &-2m(P_{x-2}\cup P_{y-3}\cup P_{c-3},k-2)
                  -tm(P_{x-1}\cup P_{y-2}\cup P_{c-2},k-1)\\
                  &-tm(P_{x-2}\cup P_{y-3}\cup P_{c-2},k-2)
                  -tm(P_{x-2}\cup P_{y-2}\cup P_{c-3},k-2)\\
                  =&2m(P_{x-3}\cup P_{y-2}\cup P_{c-2},k-1)-2m(P_{x-2}\cup P_{y-2}\cup P_{c-2},k-1)\\
                  &+m(P_{x-4}\cup P_{y-2}\cup P_{c-3},k-2)-m(P_{x-3}\cup P_{y-2}\cup P_{c-3},k-2)\\
                  &+m(P_{x-4}\cup P_{y-3}\cup P_{c-2},k-2)-m(P_{x-3}\cup P_{y-3}\cup P_{c-2},k-2)\\
                  &+2m(P_{x-3}\cup P_{y-3}\cup P_{c-3},k-2)-2m(P_{x-2}\cup P_{y-3}\cup P_{c-3},k-2)\\
                  &+tm(P_{x-2}\cup P_{y-2}\cup P_{c-2},k-1)-tm(P_{x-1}\cup P_{y-2}\cup P_{c-2},k-1)\\
                  &+tm(P_{x-3}\cup P_{y-3}\cup P_{c-2},k-2)-tm(P_{x-2}\cup P_{y-3}\cup P_{c-2},k-2)\\
                  &+tm(P_{x-3}\cup P_{y-2}\cup P_{c-3},k-2)
                  -tm(P_{x-2}\cup P_{y-2}\cup P_{c-3},k-2).
\end{align*}
By Lemmas \ref{le:matching-2} and \ref{le:matching-3}, we have $m(B_{n,x-1,y,c}^{(t+1)},k)-m(B_{n,x,y,c}^{(t)},k)<0$. So $ME(B_{n,x-1,y,c}^{(t+1)})$ $<ME(B_{n,x,y,c}^{(t)}).$
\hfill$\blacksquare$
\begin{figure}[htbp]
\begin{centering}
\includegraphics[scale=0.8]{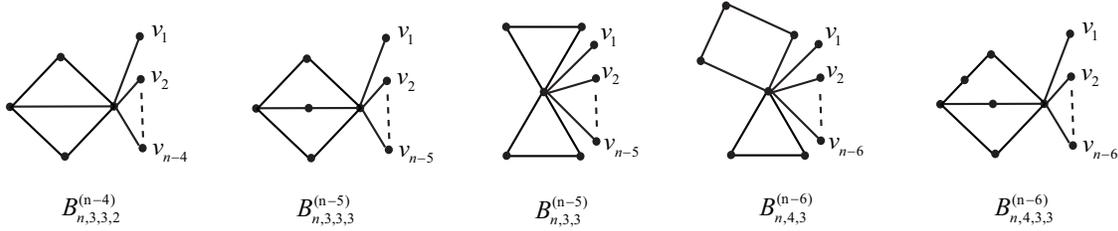}
\caption{Bicyclic graphs $B_{n,3,3,2}^{(n-4)},B_{n,3,3,3}^{(n-5)},
B_{n,3,3}^{(n-5)},B_{n,4,3}^{(n-6)} and B_{n,4,3,3}^{(n-6)}$.}\label{matching:Fig-4}
\end{centering}
\end{figure}

\noindent\begin{theorem}\label{th:matching-6}
Let $G\in B_{n}$. For $n>5$, if $G$ is not isomorphic to any member in $\{B_{n,3,3,2}^{(n-4)},$ $B_{n,3,3,3}^{(n-5)},
B_{n,3,3}^{(n-5)},B_{n,4,3}^{(n-6)}, B_{n,4,3,3}^{(n-6)}\}$; see Figure \ref{matching:Fig-4}, then we have
$$ME(G)>ME(B_{n,4,3,3}^{(n-6)})>ME(B_{n,4,3}^{(n-6)})>ME(B_{n,3,3}^{(n-5)})>ME(B_{n,3,3,3}^{(n-5)})>ME(B_{n,3,3,2}^{(n-4)}).$$
\end{theorem}
\noindent {\bf Proof.} Note that
 \begin{align*}
                  &m(B_{n,3,3,2}^{(n-4)},1)=n+1,m(B_{n,3,3,2}^{(n-4)},2)=2n-6, m(B_{n,3,3,2}^{(n-4)},k)=0\ for\  k\geq3;\\
 \end{align*}
 \begin{align*}
                  &m(B_{n,3,3,3}^{(n-5)},1)=n+1, m(B_{n,3,3,3}^{(n-5)},2)=3n-9, m(B_{n,3,3,3}^{(n-5)},k)=0\ for\  k\geq3;\\
                  &m(B_{n,3,3}^{(n-5)},1)=n+1, m(B_{n,3,3}^{(n-5)},2)=2n-5, m(B_{n,3,3}^{(n-5)},3)=n-5,  \\
                  & m(B_{n,3,3}^{(n-5)},k)=0\ for\  k\geq4;\\
                  &m(B_{n,4,3}^{(n-6)},1)=n+1, m(B_{n,4,3}^{(n-6)},2)=3n-8, m(B_{n,4,3}^{(n-6)},3)=2n-10, \\
                  & m(B_{n,4,3}^{(n-6)},k)=0\ for\  k\geq4;\\
                  &m(B_{n,4,3,3}^{(n-6)},1)=n+1, m(B_{n,4,3,3}^{(n-6)},2)=4n-13, m(B_{n,4,3,3}^{(n-6)},3)=2n-10,  \\
                  & m(B_{n,4,3,3}^{(n-6)},k)=0   \ for\  k\geq4.
\end{align*}

For $n>5$, $ m(B_{n,3,3,3}^{(n-5)},2)-m(B_{n,3,3,2}^{(n-4)},2)=n-3>0$. Thus $ME(B_{n,3,3,3}^{(n-5)})>ME(B_{n,3,3,2}^{(n-4)})$.

Similarly, we have $ME(B_{n,4,3,3}^{(n-6)})>ME(B_{n,4,3}^{(n-6)})>ME(B_{n,3,3}^{(n-5)})$.

In addition, we have
 \begin{align*}
                  &\alpha(B_{n,3,3,3}^{(n-5)},x)=x^{n}-(n+1)x^{n-2}+(3n-9)x^{n-4};\\
                  &\alpha(B_{n,3,3}^{(n-5)},x)=x^{n}-(n+1)x^{n-2}+(2n-5)x^{n-4}-(n-5)x^{n-6}.
\end{align*}

By the definition of matching energy, we have
 \begin{align*}
                  &ME(B_{n,3,3,3}^{(n-5)})=2\sqrt{\frac{n+1+\sqrt{(n+1)^{2}-4(3n-9)}}{2}}+2\sqrt{\frac{n+1-\sqrt{(n+1)^{2}-4(3n-9)}}{2}};\\
                  &ME(B_{n,3,3}^{(n-5)})=2+2\sqrt{\frac{n+\sqrt{n^{2}-4(n-5)}}{2}}+2\sqrt{\frac{n-\sqrt{n^{2}-4(n-5)}}{2}}.
\end{align*}
By a directly calculation, we have $ME(B_{n,3,3}^{(n-5)})>ME(B_{n,3,3,3}^{(n-5)})$.

Hence, $ME(B_{n,4,3,3}^{(n-6)})>ME(B_{n,4,3}^{(n-6)})>ME(B_{n,3,3}^{(n-5)})>ME(B_{n,3,3,3}^{(n-5)})>ME(B_{n,3,3,2}^{(n-4)}).$
\hfill$\blacksquare$

\end{document}